\documentclass[oneside,notitlepage,12pt]{amsart}

\usepackage{amssymb}
\usepackage[leqno]{amsmath}
\usepackage{amsfonts}
\usepackage{amsopn}
\usepackage{amstext}
\usepackage{amsthm}

\usepackage{verbatim}
\usepackage[colorlinks, backref]{hyperref}

\usepackage{calrsfs}
\usepackage{fourier-orns}
\usepackage{clock} 
\usepackage[alpine, weather]{ifsym}

\providecommand{\cal}{\mathcal}
\renewcommand{\Bbb}{\mathbb}

\newenvironment{pf}{\begin{proof}}{\end{proof}}

\newcommand{\Ef}{{\cal{F}}}

\newcommand{\Nat}{{\Bbb{N}}}

\newcommand{\eps}{\varepsilon}
\renewcommand{\phi}{\varphi}
\renewcommand{\rho}{\varrho}

\newcommand{\norm}[1]{\|#1\|}

\newcommand{\rest}{\restriction}

\newcommand{\Ntr}{n\in{\Bbb{N}}}
\newcommand{\loe}{\leq}
\newcommand{\goe}{\geq}

\newcommand{\subs}{\subseteq}
\newcommand{\sups}{\supseteq}

\renewcommand{\iff}{\Longleftrightarrow}

\newcommand{\dist}{\operatorname{dist}}

\newcommand{\G}{{\mathbb G}}

\newtheorem{tw}{Theorem}

\newtheorem{lm}{Lemma}
\newtheorem{prop}{Proposition}

\theoremstyle{definition}

\newcommand{\setof}[2]{\{#1\colon #2\}}

\newcommand{\sett}[2]{\{#1\}_{#2}}
\newcommand{\sn}[1]{\{#1\}} 
\newcommand{\map}[3]{#1\colon #2 \to #3} 
\newcommand{\img}[2]{#1[#2]} 

\newcommand{\fra}{Fra\"iss\'e}

\newcommand{\U}{\mathbb U}

\newcommand{\cmp}{\circ} 

\newcommand{\Gurarii}{Gurari\u\i}
\newcommand{\BMG}{\operatorname{BM}(\mathcal B)}
\newcommand{\BM}{\operatorname{BM}}

\title{Game-theoretic characterization of the Gurarii space}
\author{Wies{\l}aw Kubi\'s}

\address{Institute of Mathematics, Czech Academy of Sciences (CZECHIA)
}
\address{Faculty of Natural Sciences, College of Science, Cardinal Stefan Wyszy\'nski University, Warsaw (POLAND)
}

\date{\clocktime\today}

\begin{document}

\begin{abstract}
We present a simple and natural infinite game building an increasing chain of finite-dimensional Banach spaces. We show that one of the players has a strategy with the property that, no matter how the other player plays, the completion of the union of the chain is linearly isometric to the \Gurarii\ space.
\end{abstract}


\thanks{Research supported by GA\v CR grant 17-27844S and RVO 67985840 (Czechia).}

\maketitle

\section{Introduction}

We consider the following game. Namely, two players (called \emph{Eve} and \emph{Odd}) alternately choose finite-dimensional Banach spaces $E_0 \subs E_1 \subs E_2 \subs \cdots$, with no additional rules. For obvious reasons, Eve should start the game. The result is the completion of the chain $\bigcup_{\Ntr}E_n$.
We shall denote this game by $\BMG$.
This is in fact a special case of an abstract Banach-Mazur game studied recently in~\cite{KubBM}. In model theory, this is sometimes called the $\forall \exists$-game, see~\cite{Hodges-games}.
Main result:

\begin{tw}\label{ThmBMG}
	There exists a unique, up to linear isometries, separable Banach space $\G$ such that Odd has a strategy $\Sigma$ in $\BMG$ leading to $\G$, namely, the completion of every chain resulting from a play of $\BMG$ is linearly isometric to $\G$, assuming Odd uses strategy $\Sigma$, and no matter how Eve plays.
	
	Furthermore, $\G$ is the \Gurarii\ space.
\end{tw}

The result above may serve as a strong argument that the \Gurarii\ space (see the definition below) should be considered as one of the classical Banach spaces. Indeed, Theorem~\ref{ThmBMG} is completely elementary and can even be presented with no difficulties to undergraduate students who know the very basic concepts of Banach space theory.

It turns out that the \Gurarii\ space $\G$ (constructed by \Gurarii\ in 1966) is not so well-known, even to people working in functional analysis. The reason might be that this is a Banach space constructed usually by some inductive set-theoretic arguments, without providing any concrete formula for the norm. Furthermore, the fact that $\G$ is actually unique up to linear isometries was proved by Lusky~\cite{Lusky} only ten years after \Gurarii's work~\cite{Gur}.
Elementary proof of the uniqueness of $\G$ has been found recently by Solecki and the author~\cite{KubSol}. Theorem~\ref{ThmBMG} offers an alternative argument, still using the crucial lemma from~\cite{KubSol}.

In fact, uniqueness of a space $\G$ satisfying the assertion of Theorem~\ref{ThmBMG} is almost trivial: If there were two Banach spaces $G_0$, $G_1$ in Theorem~\ref{ThmBMG}, then we can play the game so that Odd uses his strategy leading to $G_1$, while after the first move Eve uses Odd's strategy leading to $G_2$. Both players win, therefore $G_1$ is linearly isometric to $G_2$.

Below, after recalling the definition of the \Gurarii\ space, we show that it indeed satisfies the assertion of Theorem~\ref{ThmBMG}.
Finally, we discuss other variants of the Banach-Mazur game, for example, playing with separable Banach spaces or with a fixed (rich enough) subclass of finite-dimensional spaces. Again, the \Gurarii\ space is the unique object for which Odd has a winning strategy.

\section{Preliminaries}

The \emph{\Gurarii\ space} is the unique separable Banach space $\G$ satisfying the following condition:
\begin{enumerate}
	\item[(G)] For every $\eps>0$, for every finite-dimensional normed spaces $A \subs B$, every isometric embedding $\map e A \G$ has an extension $\map f B \G$ that is an $\eps$-isometric embedding, namely,
	$$(1-\eps)\norm x \loe \norm{f(x)} \loe (1+\eps)\norm x$$
	for every $x \in B$.
\end{enumerate}
As we have already mentioned, this space has been found by \Gurarii~\cite{Gur} in 1966, yet its uniqueness was proved only ten years later by Lusky~\cite{Lusky} using rather advanced method of representing matrices. Elementary proof can be found in~\cite{KubSol}.
According to \cite[Thm. 2.7]{GarKub}, the \Gurarii\ space can be characterized by the following condition:
\begin{enumerate}
	\item[(H)] For every $\eps>0$, for every finite-dimensional normed spaces $A \subs B$, for every isometric embedding $\map e A \G$ there exists an isometric embedding $\map f B \G$ such that $\norm{e - f \rest A} < \eps.$
\end{enumerate}
Actually, in the proof of equivalence (G)$\iff$(H) one has to use the crucial lemma from~\cite{KubSol}:

\begin{lm}\label{LMkkrucjall}
	Let $0 < \eps < 1$ and let $\map f X Y$ be an $\eps$-isometric embedding between Banach spaces.
	Then there exists a norm on $X \oplus Y$ such that, denoting by $\map i X {X \oplus Y}$, $\map j Y {X \oplus Y}$ the canonical embeddings, it holds that
	$$\norm{j \cmp f - i} \loe \eps.$$
\end{lm}

The proof given in~\cite{KubSol} uses functionals, however there is a direct formula for the norm on $X \oplus Y$ satisfying the assertion of Lemma~\ref{LMkkrucjall}:
$$\norm{(x,y)} = \inf \setof{\norm{x_0} + \norm{y_0} + \eps \norm{x_1}}{x = x_0 + x_1, \;\; y = y_0 + f(x_1)}.$$
Easy computation showing that it works can be found in~\cite[p. 753]{CGK}.
In fact, \cite{CGK} deals with $p$-Banach spaces; $p=1$ is our case.

By a \emph{chain} of normed spaces we mean a sequence $\sett{E_n}{\Ntr}$ such that each $E_n$ is a normed space, $E_n \subs E_{n+1}$ and the norm of $E_{n+1}$ restricted to $E_n$ coincides with that of $E_n$ for every $\Ntr$.
All mappings in this note are assumed to be linear.

\section{Proof of Theorem~\ref{ThmBMG}}

Let us fix a separable Banach space $\G$ satisfying (H). We do not assume a priori that it is uniquely determined, therefore the arguments below will also show the uniqueness of $\G$.
Odd's strategy $\Sigma$ in $\BMG$ can be described as follows.

Fix a countable set $\sett{v_n}{\Ntr}$ linearly dense in $\G$.
Let $E_0$ be the first move of Eve.
Odd finds an isometric embedding $\map {f_0}{E_0}\G$ and finds $E_1 \sups E_0$ together with an isometric embedding $\map{f_1}{E_1}{\G}$ extending $f_0$ and such that $v_0 \in \img{f_1}{E_1}$.

Suppose now that $n = 2k>0$ and $E_n$ was the last move of Eve.
We assume that a linear isometric embedding $\map{f_{n-1}}{E_{n-1}}{\G}$ has been fixed.
Using (H) we choose a linear isometric embedding $\map{f_n}{E_n}{\G}$ such that $f_n \rest E_{n-1}$ is $2^{-k}$-close to $f_{n-1}$.
Extend $f_n$ to a linear isometric embedding $\map{f_{n+1}}{E_{n+1}}{\G}$ so that $E_{n+1} \sups E_n$ and $\img{f_{n+1}}{E_{n+1}}$ contains all the vectors $v_0, \dots, v_k$. The finite-dimensional space $E_{n+1}$ is Odd's move.
This finishes the description of Odd's strategy $\Sigma$.

Let $\sett{E_n}{\Ntr}$ be the chain of finite-dimensional normed spaces resulting from a fixed play, when Odd was using strategy $\Sigma$.
In particular, Odd has recorded a sequence $\sett{\map{f_n}{E_n}{\G}}{\Ntr}$ of linear isometric embeddings such that $f_{2n+1} \rest E_{2n-1}$ is $2^{-n}$-close to $f_{2n-1}$ for each $\Ntr$.
Let $E_\infty = \bigcup_{\Ntr}E_n$.
For each $x \in E_\infty$ the sequence $\sett{f_n(x)}{\Ntr}$ is Cauchy, therefore we can set $f_\infty(x) = \lim_{n \to \infty} f_n(x)$, thus defining a linear isometric embedding $\map{f_\infty}{E_\infty}{\G}$.
The assumption that $\img{f_{2n+1}}{E_{2n+1}}$ contains all the vectors $v_0,\dots, v_n$ ensures that $\img{f_\infty}{E_\infty}$ is dense in $\G$.
Finally, $f_\infty$ extends to a linear isometry from the completion of $E_\infty$ onto $\G$.
This completes the proof of Theorem~\ref{ThmBMG}.

\section{Playing with a subclass of finite-dimensional spaces}\label{Sectibwe}

It is natural to ask whether Theorem~\ref{ThmBMG} remains true when the game is restricted to a rich enough subclass of finite-dimensional normed spaces.
Of course, the minimal assumption on the class must be the existence of a chain whose completion is the \Gurarii\ space.
It turns out that this is sufficient.

Let $\Ef$ be a class of finite-dimensional normed spaces, closed under isometries. Namely, if $E \in \Ef$ and $E'$ is linearly isometric to $E$, then $E' \in \Ef$.
We say that $\Ef$ is \emph{dominating} (in the class of all finite-dimensional spaces) if for every $E \in \Ef$, for every isometric embedding $\map e E X$ with $X$ finite-dimensional, for every $\eps > 0$ there exists an $\eps$-isometric embedding $\map f X F$ such that $F \in \Ef$ and $f \cmp e$ is an isometric embedding.
Note that if $\sett{F_n}{\Ntr}$ is a chain of finite-dimensional subspaces of the \Gurarii\ space whose union is dense, then the class $\Ef$ consisting of all spaces linearly isometric to some $F_n$ is dominating.

The game $\BM(\Ef)$ is defined precisely in the same way as $\BMG$, simply restricting the class of spaces to $\Ef$.

\begin{tw}\label{ThmBMFspecjal}
	Let $\Ef$ be a dominating class of finite-dimensional normed spaces. Then Odd has a strategy $\Sigma$ in $\BM(\Ef)$ leading to the \Gurarii\ space $\G$.
	Namely, the completion of every chain resulting from a play of $\BM(\Ef)$ is linearly isometric to $\G$ when Odd uses strategy $\Sigma$.
\end{tw}

\begin{pf}
	The strategy is a suitable adaptation of the one from the proof of Theorem~\ref{ThmBMG}.
	Fix a linearly dense set $\sett{v_n}{\Ntr}$ in $\G$ such that $\norm{v_i} = 1$ for $i \in \Nat$.
	Suppose $n = 2k \goe 0$ and $E_n \in \Ef$ was the last move of Eve.
	We assume that a linear isometric embedding $\map{f_{n-1}}{E_{n-1}}{\G}$ has been defined, where $f_{-1} = 0$ and $E_{-1} = \sn0$.
	Using (H) we choose an isometric embedding $\map{f_n}{E_n}{\G}$ such that $f_n \rest E_{n-1}$ is $2^{-k}$-close to $f_{n-1}$.
	Extend $f_n$ to a linear isometric embedding $\map g X \G$ so that $X \sups E_n$ is finite-dimensional and $\{v_0, \dots, v_k\} \subs \img g X$.
	We need to ``correct" $X$ so that it becomes a member of $\Ef$.
	Using the fact that $\Ef$ is dominating, we find a $2^{-(k+1)}$-isometric embedding $\map s X F$ such that $F \in \Ef$ and $s \rest E_{n}$ is isometric. We may assume that $X \subs F$ and $s$ is the inclusion.
	We set $E_{n+1} := F$. This finishes the description of Odd's strategy, yet for the inductive arguments we still need to define the embedding $f_{n+1}$.
	
	Using Lemma~\ref{LMkkrucjall}, we find isometric embeddings $\map i X Z$, $\map j F Z$ such that $Z$ is finite-dimensional and $\norm{j \cmp s - i} \loe 2^{-(k+1)}$. Using (H), we find an isometric embedding $\map h Z \G$ such that $\norm{h \cmp i - g} \loe 2^{-(k+1)}$.
	We set $f_{n+1} := h \cmp j$. 
		
	Note that $f_{n+1} \rest X = h \cmp j \rest X = h \cmp j \cmp s$, therefore
	$$\norm{f_{n+1} \rest X - g} \loe \norm{h \cmp j \cmp s - h \cmp i} + \norm{h \cmp i - g} \loe 2^{-(k+1)} + 2^{-(k+1)} = 2^{-k}.$$
	Thus $\norm{f_{n+1} \rest E_n - f_n} \loe 2^{-k}$.
	Furthermore, if $v_i = g(x_i)$ then $\norm{f_{n+1}(x_i) - v_i} = \norm{f_{n+1}(x_i) - g(x_i)} \loe 2^{-k} \norm{x_i} = 2^{-k}$, showing that $\dist(v_i, \img{f_{n+1}}{E_{n+1}}) \loe 2^{-k}$ for $i \loe k$.
	
	Let $\sett{E_n}{\Ntr} \subs \Ef$ be the chain resulting from a play when Odd was using the strategy described above.
	In particular, we have a sequence $\sett{\map{f_n}{E_n}{\G}}{\Ntr}$ of linear isometric embeddings converging uniformly to an isometric embedding $\map{f_\infty}{E_\infty}{\G}$, where $E_\infty = \bigcup_{\Ntr}E_n$. Finally, $E_\infty$ is dense in $\G$, because
	$$\lim_{n \to \infty} \dist(v_i, \img{f_n}{E_n}) = 0$$ for each $i \in \Nat$. It follows that the unique extension of $f_\infty$ to the completion of $E_\infty$ is an isometry onto $\G$. This completes the proof.	
\end{pf}

An immediate corollary to Theorem~\ref{ThmBMFspecjal} is that if $\Ef$ is a dominating class of finite-dimensional normed spaces, then there exists a chain in $\Ef$ whose union is isometric to a dense subspace of the \Gurarii\ space.
Another corollary is the known fact that the \Gurarii\ space contains a chain of finite-dimensional $\ell_\infty$-spaces with a dense union, as the class of all such spaces is easily seen to be dominating.

\section{Final remarks}

Below we collect some comments around Theorem~\ref{ThmBMG}.

\subsection*{Universality.} It has been noticed by \Gurarii\ that $\G$ contains isometric copies of all separable Banach spaces. In fact, the space $\G$ can be constructed in such a way that it contains any prescribed separable Banach space, e.g., the space $C([0,1])$, which is well-known to be universal.
The paper~\cite{KubSol} contains a more direct and elementary proof of the isometric universality of $\G$.
The main result of this note offers yet another direct proof (cf. \cite[Thm. 10]{KubBM}).

Namely, fix a separable Banach space $X$ and fix a chain $\sett{X_n}{\Ntr} \subs X$ of finite-dimensional spaces whose union is dense in $X$. 
We describe a strategy of Eve that leads to an isometric embedding of $X$ into $\G$.
Specifically, Eve starts with $E_0 := X_0$ and records the identity embedding $\map{e_0}{X_0}{E_0}$.
Once Odd has chosen $E_n$ with $n = 2k+1$, having recorded a linear isometric embedding $\map{e_k}{X_k}{E_{n-1}}$, Eve finds $E_{n+1} \sups E_n$ such that there is a linear isometric embedding $\map{e_{k+1}}{X_{k+1}}{E_{n+1}}$ extending $e_k$. This is her response to $E_n$.
The only missing ingredient showing that such a strategy is possible is the \emph{amalgamation property} of finite-dimensional normed spaces:

\begin{lm}
	Let $\map f Z X$, $\map g Z Y$ be linear isometric embeddings of Banach spaces.
	Then there are a Banach space $W$ and linear isometric embeddings $\map{f'}{X}{W}$, $\map{g'}{Y}{W}$ such that $f' \cmp f = g' \cmp g$.
	If $X$, $Y$ are finite-dimensional then so is $W$.
\end{lm}

The above lemma belongs to the folklore and can be found in several texts, e.g., \cite{GarKub} or \cite{ACCGM}.

In any case, when Eve uses the strategy described above and Odd uses a strategy leading to the \Gurarii\ space, Eve constructs a linear isometric embedding $\map{e}{X}{\G}$ such that $e \rest X_n = e_n$ for every $\Ntr$.
This shows that $\G$ is isometrically universal in the class of all separable Banach spaces.

\subsection*{Playing with separable spaces.}
It is natural to ask what happens if both players are allowed to choose infinite-dimensional separable Banach spaces.
As it happens, in this case Odd has a very simple tactic (i.e. a strategy depending only on the last move of Eve) again leading to the \Gurarii\ space. This follows immediately from the following

\begin{prop}[{\cite[Lemma 3.3]{GarKub}}]
	Let $\sett{G_n}{\Ntr}$ be a chain of Banach spaces such that each $G_n$ is linearly isometric to the \Gurarii\ space. Then the completion of the union $\bigcup_{\Ntr}G_n$ is linearly isometric to the \Gurarii\ space.
\end{prop}

Thus, knowing that $\G$ contains isometric copies of all separable Banach spaces, Odd can always choose a space linearly isometric to $\G$, so that the resulting chain consists of \Gurarii\ spaces.

\subsection*{Other variants of the game.}
It is evident that the Banach-Mazur game considered in this note can be played with other mathematical structures. The work~\cite{KubBM} discusses this game in model theory, showing that Odd has a winning strategy leading to the so-called \emph{\fra\ limit} of a class of structures (which exists under some natural assumptions).
Another variant of this game appears when finite-dimensional normed spaces are replaced by finite metric spaces.
Almost the same arguments as in the proof of Theorem~\ref{ThmBMG} show that Odd has a strategy leading to the \emph{Urysohn space}~\cite{Urysohn}, the unique complete separable metric space $\U$ satisfying the following condition:
\begin{enumerate}
	\item[(U)] For every finite metric spaces $A \subs B$, every isometric embedding $\map{e}{A}{\U}$ can be extended to an isometric embedding $\map f B \U$.
\end{enumerate}
It turns out that $\U$ is uniquely determined by a weaker condition (analog of (H)) asserting  that $f$ is $\eps$-close to $e$ with arbitrarily small $\eps>0$, not necessarily extending $e$.
An analog of Theorem~\ref{ThmBMG} is rather obvious; the proof is practically the same as in the case of normed spaces, simply replacing all phrases ``finite-dimensional" by ``finite" and deleting all adjectives ``linear".

\subsection*{Strategies vs. tactics.}
The proof of Theorem~\ref{ThmBMG} (as well as of Theorem~\ref{ThmBMFspecjal}) actually gives a Markov strategy, that is, a strategy depending only on the step $n$ and the last move of Eve.
When playing with separable spaces, Odd has a tactic, that is, a strategy depending on the last Eve's move only (such a strategy is also called \emph{stationary}).
We do not know whether Odd has a winning tactic in the Banach-Mazur game played with finite-dimensional normed spaces or finite metric spaces, where ``winning" means obtaining the \Gurarii\ space or the Urysohn space, respectively.

\end{document}